\documentclass[12pt,a4paper]{article}
 \usepackage{epsfig,amsmath,amstext,amssymb,amscd}

\begin{document}

\title{Mathematics of Plott choice functions}
\author{    V.Danilov\thanks{The support of the RFFI 00-15-98873 grant  is grateful
acknowledged} \ and G.Koshevoy\thanks{To whom the correspondence
should to be send. The support of LIFR MIIR is grateful
acknowledged}\\ Central Institute of Economics and Mathematics\\
Nakhimovskii pr. 47, 117418 Moscow, Russia\\ e-mail: \{danilov,
koshevoy\}@cemi.rssi.ru }
\date{}

                         \maketitle

\begin{abstract}
This paper is devoted to a study of mathematical structures arising
from choice functions satisfying  the {\em path independence}
property (Plott functions). We  broaden the notion of a choice
function by allowing of empty choice. This enables us to define a
lattice structure on the set of Plott functions. Moreover, this
lattice is functorially dependent on its base. We introduce a
natural convex structure on the set of linear orders (or words) and
show that Plott functions are in one-to-one correspondence with
convex subsets in this set of linear orders. That correspondence is
compatible with both lattice structures.\medskip

\end{abstract}\medskip

                         \hfill  {\em Dedicated to the memory of Andrei Malishevski }

\section{    Introduction}
The paper is devoted to a study of choice functions satisfying the
{\em path independence} property. This property was introduced by
Plott \cite{Plott} and we will call such functions {\em Plott
functions}. He considered the concept of ``path independence'' as a
mean of weakening the condition of rationality. The issue of
''rationality of the choice'' was the main theme of investigations
of Plott functions (see, for example, \cite{A-Al,A-M,Moulin,Nehr}).
Because we do not discuss this issues. Our purpose is to show a very
nice and rich mathematics related to Plott functions.

There are two starting points for this mathematics. The first one is
the Aizerman-Malishevski theorem which states that the class of
Plott functions is identical to the class of joint-extremal choice
functions. This theorem was announced in \cite{A-M}; apparently its
first proof has appeared in \cite{Moulin}. The next important step
was made by Koshevoy \cite{Kosh}. He observed that Plott functions
can be identify with the so-called convex geometries (these
geometries were introduced in \cite{Ed}). This enables to give a
transparent explanation of the Aizerman-Malishevski theorem. More
important was that this correspondence hints on hidden mathematical
structures of Plott functions.

One circumstance impedes progress in this direction. Historically,
mainly due to psychological reasons, only choice functions with
non-empty values have been considered. The elimination of this
non-essential restriction allows us to define a lattice structure on
the set of Plott functions. To our mind, this lattice is the third
(or the fourth) most significant lattice associated with a finite
set, after the Boolean lattice of subsets and the semimodular
lattice of partitions.

The set of primitive (or join-irreducible) elements of this
lattice consists of linear Plott functions. These functions
correspond to objects which are very close to linear orders;
namely, to linear orders on subsets. To give such an order on a
set $X$ is the same as to give a simple (without repeated letters)
word over the alphabet $X$. Let ${\bf SW}(X)$ denote the set of
simple words. We define a convex structure on the set ${\bf
SW}(X)$ and show that Plott functions can be identified with
convex subsets of ${\bf SW}(X)$.

Another interest property of Plott functions is its functoriality
with respect to a base change. If $\phi:X \to Y$ is a mapping of
sets and $f$ is a Plott function on $X$ then we define a Plott
function $\phi_*(f)$ on $Y$, the direct image of $f$. The direct
image commutes with the join and pushes forward linear Plott
functions into linear Plott functions. In terms of words this
means that the corresponding map $\phi_{\sharp}:{\bf SW}(X) \to
{\bf SW}(Y)$ sends  convex sets into convex sets. Similarly one
can define the inverse image of Plott functions.

An important application of the direct image is related to the issue
of  rationalization. We show that every Plott function on a set $X$
possesses a rationalization by some partial order on some
superset $X'\to X$. Moreover, we show that there exists a canonical
(and minimal) rationalization.

The paper is organized as follows. In Section 2 we define Plott
functions, give several important examples, and introduce the
lattice of Plott functions. In Section 3 we recap the
correspondence between Plott functions and convex geometries. In
Section 4 we introduce the notion of a support of a Plott
function. The central role in our analysis of Plott functions
plays linear Plott functions (or simple words). In Section 5 we
associate to a Plott function its basement, a special subset of
linear Plott functions.  Sections 6 is devoted to the direct image
of Plott functions. In Section 7 we discuss the notion of a
superset rationalization and construct the canonical
rationalization for each Plott function. In section 8 we introduce
the shuffles and the melanges of simple words. Using the melange
operation, we introduce in Section 8 a natural convex structure on
the set ${\bf SW}(X)$ of simple words. This structure enables us
to identify Plott functions with convex subsets of ${\bf SW}(X)$.
In Section 9 we study the inverse image of Plott functions.

\section{The Plott choice functions}
     Let  $X$ be a finite set. A {\em choice function} on $X$
is a map $f\colon 2^X \to 2^X$ such that $f(A)\subset A$ for any
$A\subset X$\footnote{$A\subset B$ denotes the inclusion of sets.}.

A choice function $f$ is a {\em Plott function} (or satisfies the
{\em path independence} condition), if for any subsets $A$ and  $B$
of $X$
$$
                        f(A\cup B)=f(f(A)\cup B).
$$

Let us begin with examples of constructions of Plott
functions.\medskip

{\bf   Example 1.} Let $R$ be a partial order on  $X$ (i.e. a
reflexive, transitive, and anti-symmetric binary relation). By
$Max(A)$ (or $Max(R|A)$) we denote the set of best elements in $A$
under $R$. It is easy to check that the choice function $f_R$,
$f_R(A)=Max(R|A)$, is a Plott function.\medskip

The following instance  of this construction is important for us. It
is the case, when  $R=\le $ is a linear order, that is a complete
partial order. To define a linear order is the same as to define an
enumeration of elements of $X$, i.e  a bijection $\nu \colon
\{1,...,n\} \to X$. Here $\nu (1)$ denotes the best element of $X$,
$\nu (2)$ the second best and so on all the way to the worst element
$\nu (n)$. The choice under maximization of a linear order is a
singleton-valued choice function, i.e. $|f(A)|=1$ for any non-empty
$A\subset X$. The reverse is also true, that is a singleton-valued
Plott choice function is the maximization choice function of a
linear order.\medskip

    {\bf  Example 2.} Let $\le$ be a linear order on $X$. The
following choice function $f(A)=\{\max(A),\min(A)\}$ is a Plott
function.\medskip

More general, let  $(\le_i, i\in I)$ be a family of linear orders on
$X$. Then a  {\em joint-extremal} choice function $f$ is given by
the union of choices made under maximization for each individual
order $\le_i$. That is
$$
                        f(A)=\cup_{i\in I} \max(\le_i |A).
$$
It is easy to check that so defined function is a Plott function.
Moreover the following fact is well known \cite{Blair}:\medskip

    {\bf Lemma 1.} {\em Let  $f$ and $g$ be any two Plott functions on
$X$. Then the function $f\cup g$ given by $(f\cup g)(A)=f(A)\cup
g(A)$, $A\subset X$, is a Plott function.} $\blacksquare$\medskip

A less trivial assertion is that any Plott function with non-empty
values ($f(A)\neq\emptyset$ for any $A\neq\emptyset$) is a
joint-extremal choice function. This basic fact about Plott
functions was announced in \cite{A-M}. There are different ways to
prove it and we will discuss one of them below.
\medskip

{\bf     Example 3.}  Let $\le$ be a linear order on $X$, and our
choice from a set $A$ consists of the $k$ best elements of $A$ with
respect to this order. One can check that this is a Plott function.
It is a useful exercise to find an implementation of this choice
function as a joint-extremal choice function (even for $k=2$; see
\cite{gree} for general $k$).\medskip

   {\bf  Example 4.} Let $S$ be a subset of  $X$. Let us associate
to this subset the following choice function: ${\bf 1}_S (A)=A\cap
S$. One can check that this choice function is a Plott function
(which allows of empty choices).\medskip

The latter example has some interesting features for us. It hints us
an idea that in some sense Plott functions resemble subsets. We will
see later that a Plott function defines a certain subset of $X$, the
support of this function.  Of course, a Plott function is more
subtle object than its support (see Section 5 for more adequate
object), but this analogy with subsets suggests  to adopt a
functorial view on Plott functions.\medskip

Let us  introduce some notations. Denote by ${\bf CF}(X)$ the set
of choice functions on $X$ and by ${\bf PF}(X)$ the subset of
Plott functions. The set ${\bf CF}(X)$ (as well as ${\bf PF}(X)$)
has a natural partial order. Namely, for choice functions  $f$ and
$g$ on $X$, we write $f\le g$ if $f(A)\subset g(A)$ for any
$A\subset X$. The poset ${\bf CF}(X)$ is a lattice with $f\cup g$
and $f\cap g$ as the join and the meet of $f$ and $g$,
respectively. According to Lemma 1, in the sub-poset ${\bf PF}(X)$
any pair of elements $f$ and $g$ has the least upper bound $f\vee
g=f\cup g$. Moreover, since ${\bf PF}(X)$ has the least element
${\bf 0}$ (${\bf 0}(A)=\emptyset $ for any $A\subset X$), the
poset ${\bf PF}(X)$ is in fact a lattice. However, the greatest
lower bound $f\wedge g$ in ${\bf PF}(X)$ might differs of $f\cap
g$; it is equal to the greatest Plott function dominated by the
$f\cap g$. Later on we shall present more explicit description of
this greatest lower bound.

More generally, for any choice function $f$ there exists the maximal
Plott function $f^{\sharp}$ such that $f^{\sharp}\le f$. One can say
that $f^{\sharp}$ is the {\em Plottization} of $f$. So, $f\wedge g$
is the Plottization of $f\cap g$.

In the following figure we draw the lattice ${\bf PF}(X)$ for the
set $X$ of three elements $x,y,z$. In the centre of the picture is
the largest Plott function {\bf 1}. The smallest Plott function {\bf
0} is not drawn at all; it is situated  ``at the infinity.'' White
circles denote join-irreducible Plott functions. Words near the
circles represent the corresponding Plott functions; see Section 5.

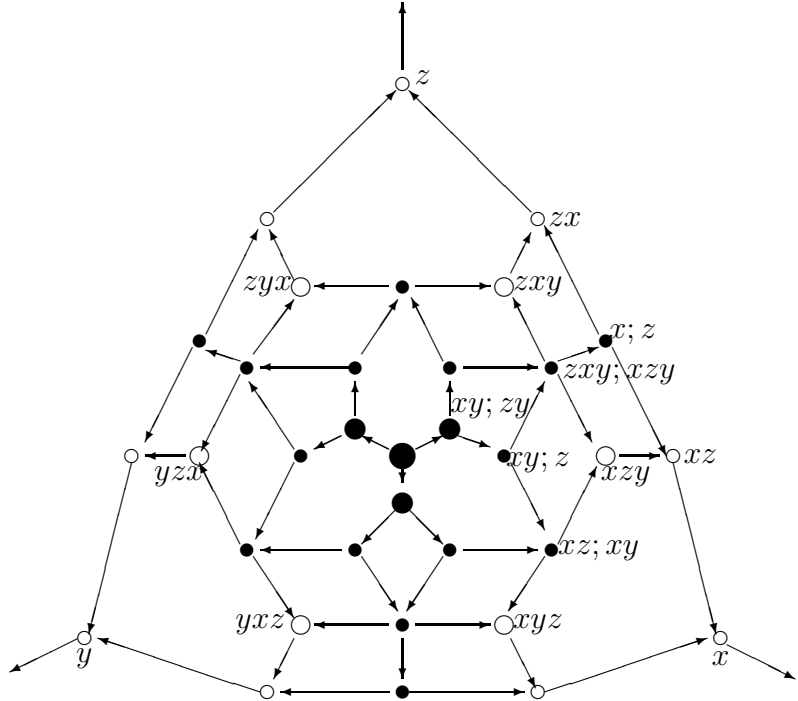
\unitlength=.900mm \special{em:linewidth 0.4pt}
\linethickness{0.4pt}
\begin{picture}(128.00,118.00)
\put(70.00,45.00){\circle*{4.00}} \put(77.00,49.00){\circle*{3.00}}
\put(63.00,49.00){\circle*{3.00}} \put(70.00,38.00){\circle*{3.00}}
\put(55.00,20.00){\circle{2.83}} \put(85.00,20.00){\circle{2.83}}
\put(100.00,45.00){\circle{2.83}} \put(85.00,70.00){\circle{2.83}}
\put(55.00,70.00){\circle{2.83}} \put(40.00,45.00){\circle{2.83}}
\put(70.00,70.00){\circle*{2.00}} \put(92.00,58.00){\circle*{2.00}}
\put(92.00,31.00){\circle*{2.00}} \put(47.00,31.00){\circle*{2.00}}
\put(47.00,58.00){\circle*{2.00}} \put(70.00,20.00){\circle*{2.00}}
\put(63.00,58.00){\circle*{2.00}} \put(77.00,58.00){\circle*{2.00}}
\put(85.00,45.00){\circle*{2.00}} \put(77.00,31.00){\circle*{2.00}}
\put(63.00,31.00){\circle*{2.00}} \put(55.00,45.00){\circle*{2.00}}
\put(70.00,10.00){\circle*{2.00}} \put(100.00,62.00){\circle*{2.00}}
\put(40.00,62.00){\circle*{2.00}} \put(50.00,10.00){\circle{2.00}}
\put(90.00,10.00){\circle{2.00}} \put(110.00,45.00){\circle{2.00}}
\put(90.00,80.00){\circle{2.00}} \put(50.00,80.00){\circle{2.00}}
\put(30.00,45.00){\circle{2.00}}
\put(72.00,46.00){\vector(2,1){4.00}}
\put(70.00,43.00){\vector(0,-1){2.00}}
\put(68.00,46.00){\vector(-2,1){4.00}}
\put(77.00,51.00){\vector(0,1){5.00}}
\put(63.00,51.00){\vector(0,1){5.00}}
\put(61.00,48.00){\vector(-2,-1){4.00}}
\put(78.00,48.00){\vector(3,-1){5.00}}
\put(71.00,37.00){\vector(1,-1){5.00}}
\put(69.00,37.00){\vector(-1,-1){5.00}}
\put(79.00,58.00){\vector(1,0){11.00}}
\put(86.00,46.00){\vector(1,2){5.00}}
\put(86.00,44.00){\vector(1,-2){5.00}}
\put(79.00,31.00){\vector(1,0){11.00}}
\put(61.00,31.00){\vector(-1,0){12.00}}
\put(54.00,44.00){\vector(-1,-2){5.67}}
\put(54.00,46.00){\vector(-2,3){6.67}}
\put(62.00,58.00){\vector(-1,0){13.00}}
\put(64.00,60.00){\vector(2,3){5.33}}
\put(76.00,59.00){\vector(-1,2){4.67}}
\put(76.00,30.00){\vector(-2,-3){5.33}}
\put(64.00,30.00){\vector(2,-3){5.33}}
\put(68.00,70.00){\vector(-1,0){11.00}}
\put(72.00,70.00){\vector(1,0){11.00}}
\put(91.00,59.00){\vector(-1,2){4.67}}
\put(93.00,57.00){\vector(1,-2){5.00}}
\put(93.00,32.00){\vector(1,2){5.67}}
\put(91.00,30.00){\vector(-2,-3){5.33}}
\put(72.00,20.00){\vector(1,0){11.00}}
\put(68.00,20.00){\vector(-1,0){11.00}}
\put(48.00,30.00){\vector(2,-3){6.00}}
\put(46.00,32.00){\vector(-1,2){6.00}}
\put(46.00,57.00){\vector(-1,-2){5.67}}
\put(48.00,60.00){\vector(3,4){6.00}}
\put(46.00,59.00){\vector(-3,1){5.00}}
\put(54.00,71.00){\vector(-1,2){3.67}}
\put(86.00,72.00){\vector(1,2){3.00}}
\put(93.00,59.00){\vector(3,1){6.00}}
\put(102.00,45.00){\vector(1,0){6.00}}
\put(86.00,18.00){\vector(1,-2){3.67}}
\put(70.00,18.00){\vector(0,-1){6.00}}
\put(54.00,18.00){\vector(-1,-2){3.00}}
\put(38.00,45.00){\vector(-1,0){6.00}}
\put(41.00,63.00){\vector(1,2){7.67}}
\put(70.00,100.00){\circle{2.00}} \put(117.00,18.00){\circle{2.00}}
\put(23.00,18.00){\circle{2.00}}

\put(39.00,61.00){\vector(-1,-2){7.00}}
\put(99.00,63.00){\vector(-1,2){8.00}}
\put(101.00,61.00){\vector(1,-2){8.00}}
\put(71.00,10.00){\vector(1,0){17.00}}
\put(68.00,10.00){\vector(-1,0){16.00}}
\put(22.00,18.00){\vector(-2,-1){10.00}}
\put(118.00,17.00){\vector(2,-1){10.00}}
\put(70.00,102.00){\vector(0,1){10.00}}
\put(23.00,15.00){\makebox(0,0)[cc]{$y$}}
\put(117.00,15.00){\makebox(0,0)[cc]{$x$}}
\put(73.00,101.00){\makebox(0,0)[cc]{$z$}}
\put(90.00,70.00){\makebox(0,0)[cc]{$zxy$}}
\put(103.00,42.00){\makebox(0,0)[cc]{$xzy$}}
\put(90.00,20.00){\makebox(0,0)[cc]{$xyz$}}
\put(49.00,20.00){\makebox(0,0)[cc]{$yxz$}}
\put(37.00,42.00){\makebox(0,0)[cc]{$yzx$}}
\put(50.00,70.00){\makebox(0,0)[cc]{$zyx$}}
\put(83.00,52.00){\makebox(0,0)[cc]{$xy;zy$}}
\put(90.00,44.00){\makebox(0,0)[cc]{$xy;z$}}
\put(102.00,57.00){\makebox(0,0)[cc]{$zxy;xzy$}}
\put(99.00,31.00){\makebox(0,0)[cc]{$xz;xy$}}
\put(104.00,63.00){\makebox(0,0)[cc]{$x;z$}}
\put(94.00,80.00){\makebox(0,0)[cc]{$zx$}}
\put(114.00,45.00){\makebox(0,0)[cc]{$xz$}}
\put(30.00,44.00){\vector(-1,-4){6.33}}
\put(49.00,10.00){\vector(-3,1){24.00}}
\put(91.00,10.00){\vector(3,1){24.00}}
\put(110.00,44.00){\vector(1,-4){6.33}}
\put(89.00,81.00){\vector(-1,1){18.00}}
\put(51.00,81.00){\vector(1,1){18.00}}
\end{picture}

\section{ Connection with convex geometries}
We have seen that the set ${\bf PF}(X)$ of Plott functions on $X$ is
a lattice. However, even on the individual level, each Plott
function is canonically related to some special lattice, the
so-called convex geometry. This construction was introduced in
\cite{Kosh} (see also \cite{MojRad}). Let us recap some details of
this construction, mainly to agree on notations.

Recall, that a collection  $\mathcal F $ of subsets of $X$ (we
will refer to the elements of $\mathcal F $ as to closed or convex
subsets of $X$) forms a {\em convex geometry}, if  $\mathcal F $
is stable with respect to the intersection, contains  $X$, and
possesses the {\em Minkowski-Krein-Mil'man (MKM) property}. To
formulate the latter property we have to introduce the notion of
an extreme point. A point $a$ of a set $A\subset X$ is said to be
{\em extreme} (with respect to $\mathcal F $), if $a$ does not
belong to the closure of the set $A\setminus a$ (the closure of a
set $B\subset X$  is the minimal set of $\mathcal F$ which
contains $B$). The set of extreme points of $A$ is denoted by
$ext(A)$ (or $ext_{\mathcal F}(A)$, if we want to emphasize the
closure system we are interested in). A collection $\mathcal F$ of
subsets of $X$ is said to satisfy the MKM property if the closure
of $A$ coincides with the closure of $ext(A)$ for any $A\subset
X$. Roughly speaking, each subset has ``enough'' extreme points.
Denote by ${\bf CG}(X)$ the set of all convex geometries on $X$.

There are the following mutually inverse mappings between the
introduced sets \cite{Kosh} (the lattice isomorphisms in fact):
$$
                   {\bf PF}(X) \to {\bf CG}(X) \to {\bf
PF}(X).
$$
The mapping from Plott functions to geometries is defined by the
following rule: the closure of a set  $A$ consists of all elements
$x \in X$, such that $f(A)=f(A\cup x)$. The inverse mapping sends a
convex geometry $\mathcal F$ to the choice function $ext_{\mathcal
F}$.

This construction enables us to switch from one language to
another. Having formulated something in one language, we will
sometimes translate this into the other language, and sometimes we
will leave such a translation to the reader. For example, above we
considered the join of a pair of Plott functions. In the language
of convex geometries, the join $\mathcal F \vee \mathcal G $
consists of the sets  $F\cap G$ with $F \in \mathcal F $ and $G\in
\mathcal G $. One more example: a Plott function $f$ takes only
non-empty values if and only if $\emptyset $ is a closed set of
the corresponding convex geometry $\mathcal F $.

\section{ The support }
With each Plott function $f$ we associate a set, the {\em support}
of $f$. This set is given by the following formula:
$$
                     S= supp(f)=\{x\in X, f(x)=x\}.
                   $$

  {\bf    Proposition 1.} {\em  For any $A\subset X$, there holds
$f(A)=f(A\cap S)$.}\medskip

{\bf Proof}. We proceed by induction on the number of elements of
$A\setminus S$. If there is no such elements, then  $A\cap S=A$
and the proposition is trivially true. Let $x \notin S$ and
$A'=A\cup x$. Since $A'\cap S=A\cap S$, we have to show $f(A\cup
x)=f(A)$. According to the path independence property, we have
$$
                   f(x\cup A)=f(f(x)\cup A)=f(\emptyset \cup
A)=f(A). \qquad\qquad\qquad               \blacksquare
$$

Thus, we see that elements outside the support $S$ are irrelevant to
the choice; in terms of convex geometries $X\setminus S$ is the
minimal convex set (the closure of $\emptyset $). Actually such a
function is defined on $S$. We will develop a formalism which
provides this intuition.

It is clear that the mapping
$$
                         supp \colon {\bf PF}(X) \to 2^X
                        $$
sending a choice function to its support is a $\vee $-morphism of
lattices. That is a monotone mapping (with respect to the partial
orders on these lattices) commuting with the join operation  $\vee
$. In general case, this map does not commute with $\wedge $. There
is a natural ''reverse'' mapping (see Example 4), which also
commutes with  $\vee $.

The support of a Plott function is a simple, but a rather rough
characteristic. It is exactly in this sense that we pointed out
before that Plott functions remind us of subsets. For subsets we
might speak about their image and inverse images under mappings. In
Sections 6 and 9 we shall define direct and inverse images of Plott
functions. That means that ${\bf PF}(X)$ is functorial with respect
to the base $X$.

 \section{ Linear Plott functions and simple words}
{\bf Definition.} A Plott function  $f$ is said to be {\em linear}
if $f(A)$ contains at most one element for any $A \subset
X$.\medskip

Example 1 show how one can construct linear Plott functions by the
means of linear orders. Slightly more general construction is as
follows. Fix a subset $S\subset X$ and a linear order $\le_S$ on $S$
and define a choice function $l$ by the following formula:  for
$A\subset X$ \ $l(A)$ is the best element (with respect to $\le_S$)
of $A\cap S$ (if $A\cap S=\emptyset$ then $l(A)=\emptyset$).
Obviously, $l$ is a linear Plott function with the support $S$.

The reverse is also true. For functions with the full support $X$
this was proven by Plott \cite{Plott}. In the general case, we have
to restrict a linear Plott function to its support.

To give a linear order on a subset of $X$ is the same as to give a
simple word over the alphabet $X$. A {\em word} is a sequence
$x(1)x(2)...x(k)$ of elements $x(i)\in X$, $k=0,1,...$. A word is
{\em simple} if no letter is repeated. A simple word
$w=x(1)x(2)...x(k)$ define a linear order
$\le_w=(x(1)>x(2)>...>x(k))$ on the subset
$supp(w)=\{x(1),x(2),...,x(k)\}\subset X$. The set of simple words
is denoted by ${\bf SW}(X)$. We have seen that the set ${\bf
SW}(X)$ can be identify with the set ${\bf LPF}(X)$ of linear
Plott functions.

It is easy to describe a "linear" convex geometry corresponding to a
simple word $w=x(1)x(2)...x(k)$. It consists of subsets
$$
                      X=\mathcal L
(0)\supset \mathcal L (1)\supset ...\supset \mathcal L (k),
                    $$
where $\mathcal L (j)=X\setminus \{x(1),...,x(j)\}$, $j\le k$. Vice
versa, a chain $\mathcal L$:
$$
                      X=\mathcal L
(0)\supset \mathcal L (1)\supset ...\supset \mathcal L (k)
                    $$
of subsets is a convex geometry if and only if any two consecutive
members of this chain differ exactly by a single element.

The identification of ${\bf SW}(X)$ and ${\bf LPF}(X)$ is
compatible with orders on these sets if we set that a word is
larger  any of  its prefixes. Picture 2 illustrates the poset
${\bf SW}(X)$ for a three-elements set $X$.

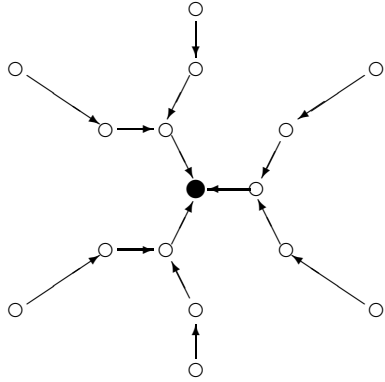
\unitlength=.8mm \special{em:linewidth 0.4pt} \linethickness{0.4pt}
\begin{picture}(100.00,73.00)(-10,0)
\put(60.00,5.00){\circle{2.00}} \put(60.00,15.00){\circle{2.00}}
\put(90.00,15.00){\circle{2.00}} \put(30.00,15.00){\circle{2.00}}
\put(45.00,25.00){\circle{2.00}} \put(55.00,25.00){\circle{2.00}}
\put(75.00,25.00){\circle{2.00}} \put(70.00,35.00){\circle{2.00}}
\put(60.00,35.00){\circle*{3.00}} \put(75.00,45.00){\circle{2.00}}
\put(55.00,45.00){\circle{2.00}} \put(45.00,45.00){\circle{2.00}}
\put(30.00,55.00){\circle{2.00}} \put(60.00,55.00){\circle{2.00}}
\put(90.00,55.00){\circle{2.00}} \put(60.00,65.00){\circle{2.00}}
\put(56.00,44.00){\vector(1,-2){3.50}}

\put(56.00,26.00){\vector(1,2){3.5}}
\put(69.00,35.00){\vector(-1,0){7.00}}
\put(60.00,63.00){\vector(0,-1){6.00}}
\put(59.00,54.00){\vector(-1,-2){3.67}}
\put(32.00,54.00){\vector(3,-2){12.00}}
\put(47.00,45.00){\vector(1,0){6.00}}
\put(88.00,54.00){\vector(-3,-2){11.00}}
\put(74.00,43.00){\vector(-1,-2){3.00}}
\put(88.00,16.00){\vector(-3,2){12.00}}
\put(74.00,26.00){\vector(-1,2){3.67}}
\put(60.00,7.00){\vector(0,1){5.50}}
\put(59.00,17.00){\vector(-1,2){3.00}}
\put(32.00,16.00){\vector(3,2){12.00}}
\put(47.00,25.00){\vector(1,0){6.00}}
\end{picture}

{\bf Definition.} The {\em basement} of a choice function $f$ is the
set $Bas(f)$ of all linear Plott functions $l$ such that $l\le
f$.\medskip

Using the identification of linear Plott functions with simple
words we can consider any basement as a subset of ${\bf SW}(X)$.
As we shall see, any Plott function $f$ is equal to the join of
all linear Plott functions of the basement of $f$. However, in the
beginning we consider several examples.\medskip

   {\bf  Example $4'$.} Let $f={\bf 1}_X$ be the identity choice
function on $X$ ($f(A)\equiv A$). Then the basement of this
function is the whole set  ${\bf LPF}(X)$ (or ${\bf SW}(X)$).
\medskip

    {\bf  Example 5.} Suppose $l_w$ is a linear Plott function
corresponding to a simple word $w=x(1)x(2)...x(k)$. The basement
of $l_w$ consists of all prefixes of $w$. \medskip

    {\bf Example $2'$.} Suppose that $X=\{a,b,c\}$ and that $f$ is the
join of two linear Plott functions corresponding to the words $abc$
and $cba$. Then the basement of $f$ consists of the following words:

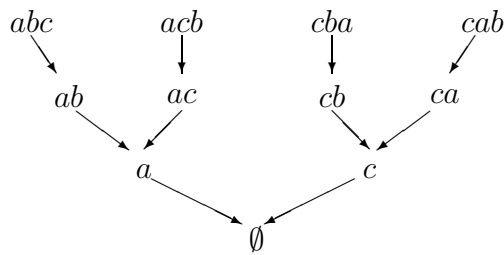
\unitlength=1mm \special{em:linewidth 0.4pt} \linethickness{0.4pt}
\begin{picture}(90.00,40.00)
\put(60.00,6.00){\makebox(0,0)[cc]{$\emptyset$}}
\put(45.00,15.00){\makebox(0,0)[cc]{$a$}}
\put(75.00,15.00){\makebox(0,0)[cc]{$c$}}
\put(35.00,25.00){\makebox(0,0)[cc]{$ab$}}
\put(50.00,25.00){\makebox(0,0)[cc]{$ac$}}
\put(70.00,25.00){\makebox(0,0)[cc]{$cb$}}
\put(85.00,25.00){\makebox(0,0)[cc]{$ca$}}
\put(30.00,35.00){\makebox(0,0)[cc]{$abc$}}
\put(50.00,35.00){\makebox(0,0)[cc]{$acb$}}
\put(70.00,35.00){\makebox(0,0)[cc]{$cba$}}
\put(90.00,35.00){\makebox(0,0)[cc]{$cab$}}
\put(30.00,33.00){\vector(2,-3){3.33}}
\put(50.00,33.00){\vector(0,-1){5.00}}
\put(70.00,33.00){\vector(0,-1){5.00}}
\put(89.00,33.00){\vector(-2,-3){3.33}}
\put(83.00,23.00){\vector(-4,-3){7.00}}
\put(70.00,23.00){\vector(1,-1){5.00}}
\put(50.00,23.00){\vector(-1,-1){5.00}}
\put(36.00,23.00){\vector(4,-3){7.00}}
\put(46.00,14.00){\vector(2,-1){12.00}}
\put(73.00,14.00){\vector(-2,-1){12.00}}
\end{picture}

\noindent where $\emptyset$ denotes the empty word.\medskip

Sometimes it is convenient to deal with special subset of
$Bas(f)$.\medskip

{\bf Definition.} The {\em socle} of a choice function $f$ is the
set $Soc(f)$ of maximal elements of $Bas(f)$.\medskip

{\bf Example $1'$.} Let $R$ be a partial order on $X$ and let $f_R$
be the Plott function arising from maximization $R$. It is easy to
check that the socle of $f_R$ consists of all linear Plott functions
corresponding to all linear extensions of $R$.\medskip

For any choice function $f$ we have $f\ge \vee(Bas(f))$. The
following result states that the equality holds for any Plott
fumction $f$\footnote{ Here is a straightforward analogy with the
standard convexity. Linear orders correspond to linear functionals
and a chain of orders corresponds to a chain of half-spaces. Similar
to a polytope which is given by the intersection of half-spaces, our
''convex'' sets will be given by the intersection of
chains-''half-spaces''. Theorem 1 is similar to the assertion that a
convex functions $f$ is supremum of linear (affine) functions which
are less than or equal to $f$.}. It generalizes and refines the
Aizerman-Malishevki theorem.\medskip

{\bf  Theorem  1.} {\em Any Plott function $f$ is equal to the join
of all linear Plott functions from $Bas(f)$. That is}
$$
f=\vee(Bas(f)).
$$

In ''geometric'' language this theorem asserts that any convex
geometry is equal to the join of its ''linear'' convex
sub-geometries. We shall prove the theorem in this geometric form.
\medskip

{\bf Proof}. Consider a convex geometry $\mathcal F $ on $X$, and
consider a maximal  chain  $\mathcal L =( L (0)\supset  L (1)\supset
...$), where all $L (k)\in \mathcal F $. We claim that $\mathcal L $
is a ''linear'' convex sub-geometry. Since it is a chain, the
linearity follows from the construction. We have to check that
adjacent members of this chain differ exactly by a single element.

So, let  $F= L (k)$ and  $G= L (k+1)$ be any two consecutive entries
of the chain $\mathcal L $. We claim that there exists an extreme
point  $x$ of the set  $F$, which does not belong to $G$. Indeed,
the set $ext(F)$ can not be a subset of $G$, otherwise by the MKM
property, this would imply that $F$ (as the closure of $ext(F)$) is
a subset of  $G$. Consider the set $F'=F\setminus x$. Since $F$ is
closed and $x$ is an extreme point of $F$, the set $F'$ is closed,
and $F'$ is located between  $F$ and $G$. Since the chain is
maximal, we have  $F'=G$, that proves our assertion.

Now the theorem is obvious. We have proven an even more general
result: a convex geometry $\mathcal F $ is equal to the union of
its ''linear'' sub-geometries. This holds because any element of
$\mathcal F $ might be included in some maximal  chain. In fact,
any chain in $\mathcal F $ might be condensed until it becomes
maximal. $\blacksquare$\medskip

   As a consequence of this theorem and Example 5 we obtain that
linear Plott functions are indeed join-irreducible elements of the
lattice ${\bf PF}(X)$.

Another consequence is that the Plottization of $f$ is equal to the join
of function from the basement of $f$.

Theorem 1 demonstrates that any Plott function can be identified
with the subset $Bas(f)\subset {\bf SW}(X)$. Later on we introduce
a convex structure on ${\bf SW}(X)$ and show that the basements of
Plott functions are precisely the convex subsets of ${\bf SW}(X)$.
For this and other aims we  study the behavior of Plott functions
with respect to exchange of the base set $X$.

\section{ The direct image}
Until now we dealt with choice functions on a fixed set $X$. Now we
compare these functions when $X$ varies.

Suppose we are given a mapping $\phi\colon   X\to Y$, and let $f$ be
a choice function $f$ on $X$. Then we may ``push forward'' $f$ from
$X$ to $Y$ (similar to subsets and linear orders). Specifically,
define a choice function $\phi_*(f)$ on $Y$, sending $B\subset Y$ to
$$
\phi_*(f)(B)=\phi (f(\phi^{-1} B)).
$$
The choice function $\phi_*(f)$ is called the {\em direct image} of
$f$ under the map $\phi$.

For each element $y\in Y$, we will refer to the elements of the
pre-image $\phi^{-1}(y)$ as to clones of $y$ and denote them
$y',y'',\ldots$ (however, some $y\in Y$ might have no clones at
all). Then an ''original element'' $y\in B$ belongs to the choice
$\phi_*(f)(B)$ if and only if at least one of the clones of $y$
belongs to the choice under $f$, i.e. $\phi^{-1}(y)\cap
f(\phi^{-1}(B))\neq\emptyset$. It is clear that the map
$\phi_*\colon {\bf CF}(X)\to {\bf CF}(Y)$ is monotone and commutes
with $\vee $.\medskip

    {\bf Example  $4''$.} Let  $f={\bf 1}_S$ be the identity
choice function with the support $S\subset X$. Then $\phi_*(f)={\bf
1}_{\phi (S)}$.\medskip

    Moreover, the support of the choice function $\phi_*(f)$ is
equal to the image of the support of $f$.\medskip

We assert that the direct image $\phi_*(l)$ of a linear Plott
function $l$ on $X$ is a linear Plott function on $Y$. Let a simple
word $w$ generate $l$. Then we explicitly write a simple word
generates $\phi_*(l)$.

For this we define a map
$$
{\bf SW}(\phi)=\psi_{\sharp}:{\bf SW}(X) \to {\bf SW}(Y).
$$
For a word $w=x(1)...x(k)$ let $\phi(w)$ denote a word
$\phi(x(1))...\phi(x(k))$. For a simple word $w$ we define
$\phi_{\sharp}(w)$ as a simplification of the word $\phi(w)$. Specifically, if $w=vx$ then
$$
\phi_{\sharp}(w)=\left\{\begin{array}{cl}
  \phi_{\sharp}(v)\phi (x) & \text{ if } \phi(x) \text{ does not occur in }
  \phi(v), \\
  \phi_{\sharp}(v) & \text{ otherwise.}
\end{array}\right.
$$

Suppose, for example, $Y=\{a,b, c, d\}$, $X=\{a',a'',b',c',c'',c'''\}$ and $\phi(a')=\phi
(a'')=a$ and so on. Then for the word $w=c''b'c'a''a'$ its direct
image is the word $\phi_* (w)=cba$.

Obviously there holds\medskip

    {\bf  Lemma 2.} {\em  Let $l_w$ be a linear Plott function
on $X$, corresponding to a simple word  $w$. Then $\phi_* (l_w)$ is
a linear Plott function on $Y$, corresponding to the word
$\phi_{\sharp} (w)$. }$\blacksquare$\medskip

In other words, we have a commutative diagram
$$
\begin{CD}
{\bf LPF}(X) @>{\phi_*}>> {\bf LPF}(X) \\
@VV{\simeq}V @VV{\simeq}V \\
{\bf SW}(X) @>{\phi_{\sharp}}>> {\bf SW}(Y)
\end{CD}
$$
where the vertical maps identify linear Plott functions and simple
words.\medskip

{\bf    Corollary.} {\em  The direct image of a Plott function is a
Plott function.}\medskip

{\bf Proof}. Because of Theorem 1, a Plott function $f$ on $X$ is
equal to the join of a set $B$ of linear Plott functions. Since
$\phi_*$ commutes with $\vee $, $\phi_*(f)$ is equal to the join
($\vee$) of the linear Plott functions $\phi_*(l)$, $l\in B$. \hfill
 $\blacksquare$\medskip

It is clear that if  $\phi   \colon     X \to Y$ and $\psi  \colon Y
\to Z$ are two maps, then there holds $(\psi \circ \phi )_* =\psi_*
\circ \phi_*$ . This implies that  {\bf PF} is a covariant functor
from the category of (finite) sets to the category of partially
ordered sets (even more, $\vee$-lattices).

     In particular, if  $X$ is a subset of  $Y$, we may speak
about a ``trivial'' extension of a choice function $f$ from $X$ to
$Y$. We denote it by  $f_Y$, and it is set by the rule
$$
                           f_Y(B)=f(B\cap X).
                        $$
Obviously, the support of $f_Y$ is a subset of $X$. We can identify
Plott functions on $X$ and Plott functions on $Y$ with the support
in $X$.\medskip

The following Theorem refines the previous corollary.\medskip

{\bf Theorem 2.} $Bas(\phi_*(f))=\phi_*(Bas(f))$.\medskip

Or, equivalently,  if we understand the basements as subsets of
${\bf SW}(X)$ and ${\bf SW}(Y)$, then
$Bas(\phi_*(f))=\phi_{\sharp}(Bas(f))$.\medskip

A proof of this theorem is proceeded in the language of convex geometrices,
and, therefore, in the beginning  we give a ''geometric'' description of
the direct image. Let $ \phi :X \to Y$ be a mapping of sets, and $
\mathcal F$ be a convex geometry on $X$. Let us define
$$
\mathcal F_Y=\{ \phi_+(F), \ F\in \mathcal F\}.
$$
Here $\phi_+$ denotes the full image, that is for a subset $A\subset
X$
$$
\phi_+(A)=\{y\in Y, \phi^{-1}(y)\subset A\}.
$$
In particular, if $y\notin \phi (X)$ then $y\in \phi_+(A)$ for any
$A \subset X$. The operation $\phi_+$ commutes with the intersection
$\cap$.\medskip

     {\bf Proposition 2.} {\em Let $\mathcal F$ be the convex
geometry corresponding to a Plott function $f$ on $X$. Then $\mathcal
F_Y$ is the convex geometry on $Y$ which corresponds to
$\phi_*(f)$.}\medskip

{\bf Proof}. A) $\mathcal{F}_Y$ is a lattice. Indeed,
$Y=\phi_+(X)$ belongs to $\mathcal F_Y$. Since $\phi_+$ commutes
with the intersection, $\mathcal F_Y$ is stable with respect to the
intersection.

B) Let us check that $\mathcal F_Y$ is a convex geometry.

     For this we have to show that if $G\supset G'$ are adjacent
closed sets of $\mathcal F_Y$ (i.e. there are no closed sets
located between  $G$ and $G'$), then they differ by a single element.
By the definition of $\mathcal F_Y$ there are convex sets $F$ and $F'$ (of $\mathcal F$)
such that $G=\phi_+(F)$ and $G'=\phi_+(F')$. Set $ F''=F\cap F'$. It
is clear that $F''\subset F$. On the other hand
$$
                   \phi_+(F'')=\phi_+(F)\cap \phi_+
(F')=G\cap G'=G'.
                 $$
Replacing, if needed, $F'$ by $F''$ we can suppose that $F\supset F'$. Let now
$F=F_0\supset ...\supset F_k=F'$ be a condensend chain of convex
sets between $F$ and $F'$. Since $G$ and $G'$ are adjacent, any set
of the form $\phi_+(F_j)$ is either $G$ or $G'$. Since $\mathcal F$
is a convex geometry then neighbor members of the chain $F_j$ differ
by a single element. Therefore $\phi_+(F_j)$ and $\phi_+(F_{j+1})$
differ by at most one element. It is now obvious that $G$ and $G'$
differ by a single element.

C)      Now we describe extreme points of an arbitrary subset
$B\subset Y$. By the definition, a point $b\in B$ is extreme if
there exists a convex set $G\subset Y$ which contains $B\setminus b$ but does
non contain $b$. That is there exists a convex set $F\subset X$ which contains $
\phi^{-1}(B\setminus b)$ but not $\phi^{-1}(b)$. That means
exactly that the fiber $\phi^{-1}(b)$ contains an extreme point of
the set $\phi^{-1}(B)$. In other words, we have proven the following
formula
$$
    ext_{\mathcal F_Y}(B)=\phi (ext_{\mathcal F}(\phi^{-1} (B)).\blacksquare
$$

{\bf Proof of Theorem 2}. The inclusion $ \supset$ is obvious. Inversely,
let $l$ be a linear Plott function on $Y$ such that $l \le
\phi_*(f)$. Consider the corresponding chain $\mathcal L=(
Y=G_0\supset G_1\supset ...\supset G_k)$ of convex sets in $Y$.
Accordingly to the step B) of the previous proposition proof, this chain can be
lifted to a chain $X=F_0\supset F_1\supset ...\supset F_k$ of convex
sets in $X$. If we condense this chain we obtain a linear
subgeometry of $\mathcal F$ which projects in $\mathcal L$.
$\blacksquare $

\section{Superset rationalizations of Plott functions }

A Plott function $f$ on  $X$ is called {\em rationalizable} by a
partial order $R$ on $X$ if, for any $A\subset X$, $f(A)$ consists of
maximal elements in $A$ with
respect to this order $R$ (see Example 1). The partial order $R$ (if
it exists) is uniquely defined by $f$ and is called the {\em
rationalization} of $f$.

Not every Plott function is rationalizable. For example, the
function from Example 2' is not rationalizable. In literature
(especially in papers by Malishevski and by Nehring), it was considered
another form of  rationalization of Plott functions using transitive
hyper-relations on $X$ (that are special relations between
$2^X$ and $X$). Here we propose a rationalization of Plott functions by transitive
relations on a "superset" $X' \to X$.\medskip

{\bf Definition}. A triple $(Y,\le ,\psi )$ (where $\le $ is a
partial order on $Y$ and $\psi  \colon   Y \to X$ is a mapping) is
said to be a {\em SS-rationalization} of $f$ if $f=\psi_*(g)$, where
$g$ is the  Plott function on $Y$ rationalizable by $\le$.\medskip

  {\bf  Example $2''$}. Let  $X=\{a,b,c\}$ and let a choice function
$f$ be the join of two linear orders, which correspond to the words
$abc$ and $cba$ (see Example  2'). Then the poset
$$
                          a' \to b' \to c'
                    $$
$$
                         c'' \to b'' \to a''
                     $$
is a rationalization of $f$ ($\psi(x')=\psi(x'')=a$,
$x\in\{a,b,c\}$). However there exists a  rationalization
of $f$ by the smaller poset
$$
                       a' \to b', c' \to b''.
                   $$

Let us explain a meaning of the SS-rationalization. A main idea of the
binary rationalization consists in that that a rationalizable choice
can be made on a pairwise base and the binary relation defined by this pairwise
choice is an order. Then an element $x$ is rejected in
$A$, if the considered opportunity set $A$ contains an element $x'$
which dominates $x$ (with respect to the order). We follow this idea
when we deal with a superset. Namely, each original
$x$ is split into several splinters, and all these splinters form a superset $X'$.
An order on $X'$ is defined such that its direct image
(under the natural mappins sending splinters to their original)
is exactly the choice function under consideration.  Now, the original $x$ is
rejected if and only if EVERY of its splinters is dominated.

In other words, each alternative has many aspects or qualities. And
if an alternative is not dominated with respect to some of its
quality, then  we have a rationale to include this alternative in
our choice. And conversely, if our alternative $x$ is dominated with
respect to any of its qualities, we have a rational to exclude it
from our choice. Of course, to exclude an alternative we have to use
several other alternatives. So, in Example $2''$ we exclude $b$ from
the set $\{a,b,c\}$ because all of its splitters $b'$ and $b''$ are
dominated.

It is easy to prove (using Theorem 1) that any Plott function
possesses (a lot of) SS-rationalizations. Conversely, an
SS-rationalizable choice function is a Plott function indeed. Thus, Plott
functions are precisely SS-rationalizable choice functions; this
assertion was proven by \cite{Moulin}. We shall not digress to prove
this assertion because we provide any Plott function with a canonical (and minimal, in a
sense) SS-rationalization. For this we define a canonical splitting of
elements of $X$ into aspects or pieces.

The following notion will be of use in this section. Let $f$ be a
Plott function on  $X$ and let $\mathcal F $ be the corresponding
convex geometry on  $X$. Recall, that elements of $\mathcal F $ we
call convex sets.\medskip

{\bf    Definition.} A closed  set  $P\in\mathcal F$ is said to be a
{\em piece} of an element  $x\in X$ if  $x\notin P$ and $x\in Q$ for
any proper closed superset  $Q\in\mathcal F$ of $P$ ($P\subset Q$
and $P\ne Q$).\medskip

It immediately follows from the definition that a piece of $x$ is a
maximal convex set (under inclusion) not containing $x$. In
particular, if $x$ does not belong to the closure (in $\mathcal F$)
of a set $A\subset X$, then there exists a piece of $x$ which
contains $A$.

A convex set  $P$ is said to be  a {\em piece} if it is a piece of
an element. It is important that any piece is exactly a piece of a
unique element of $X$.\medskip

    {\bf  Lemma 3}. {\em  Let  $P$ be a piece of  $x$. Then the set
$P\cup x$ is convex}.\medskip

{\bf Proof}. Let  $Q$ covers  $P$ in  the lattice $\mathcal F $
(i.e. for any $Z\in \mathcal F$ such that $Q\supset Z\supset P$
either $Z=Q$ or $Z=P$). Then $x\in Q$. Because $Q$ and $P$ differ by
a single element (the MKM property) , we have $Q=P\cup x$.
$\blacksquare$\medskip

{\bf     Corollary}. {\em  Let $P$ be a piece of $x$ and a piece of $x'$.
Then $x=x'$.}\medskip

In fact, the set  $P\cup x'$ is convex and contains  $P$. Hence
$P\cup x'$ contains  $x$. That implies  $x=x'$.
$\blacksquare$\medskip

Let $f$ be a Plott function. Denote by  $\mathcal P (f)$ the set of
all pieces of the corresponding convex geometry $\mathcal F$. It is
easy to see that any convex set can be given as the intersection of
some pieces. On the other hand, any piece $P$ cannot be presented as
the intersection of  convex sets, each of which properly contains $P$.
In fact, any large convex
set will contain $P\cup x$, if $P$ is a piece of $x$. Thus $\mathcal
P(f)$ is a minimal basis (with respect to the meet operation) of the
lattice $\mathcal F$.

Since each piece is a piece of a unique element, the natural map
$\phi \colon   \mathcal P (f) \to X$, which sends $P\in \mathcal P
(f)$ to the unique element $\phi (P)\in X$ for which $P$ is a piece,
is correctly defined. The inverse image  $\phi^{-1} (x)$ is
constituted of all pieces of  $x$. Because each element $x\in supp
(f)$ has a piece, we have $\phi(\mathcal P (f))=supp (f)$. \medskip

    {\bf  Example $2'''$}. Consider the choice function from Example
$2''$. The element  $a$ has exactly one piece $a'=\{b,c\}$, the
element $c$ also has a unique piece  $c'=\{a,b\}$. While the element
$b$ has two pieces  $b'=\{a\}$ and  $b''=\{c\}$.\medskip

The set  $\mathcal P (f) $ endowed with the inclusion order $\subset
$ is a poset. We state that the poset  $\mathcal P (f)$ with the
natural map $\phi  \colon   \mathcal P \to X$ is an
SS-rationalization of $f$.\medskip

    {\bf  Theorem 3}. {\em Let $f$ be a Plott function. Then the triple
$(\mathcal P (f),\subset ,\phi )$ is an SS-rationalization of $f$.
}\medskip

{\bf Proof}. Denote by $\tilde f$ the choice function on $\mathcal
P (f)$, which is rationalizable by the order $\subset$.
We have to check that $f=\phi_*(\tilde f)$,
i.e. for any subset $A\subset X$ the following equality holds
$$
                       f(A)=\phi (\tilde f(\phi^{-1}(A)).
$$
Let  $a\in f(A)$, that is $a$ is an extreme point of  $A$. By the
definition of an extreme point  $a$ does not belong to the closure
of  the set $A\setminus a$. Denote by  $P$ (any) piece of $a$, which
contains  $A\setminus a$. We claim that $P$ is not contained in any
element (of $\mathcal P (f)$) of $\phi^{-1}(A)$. In fact, assume
there exists a piece $P'$ of some element $a'$ of $A\setminus a$,
such that  $P'\supset P$. Then $a'\in A\setminus a\subset P\subset
P'$, that contradicts to $a'\notin P'$. Thus $a\in (\phi_*(\tilde
f))(A)$.

Now, let  $a\notin f(A)$. Let us show that any piece  $P$ of  $a$ is
a subset of some element of $\phi^{-1}(A\setminus a)$. Because $f$
is a Plott function, $a$ does not belong to the set of extreme
points of $A$. Therefore the closure of $A\setminus a$ contains $a$.
Thus $A\setminus a$ is not a subset of  $P$ (otherwise, the closure
of $A\setminus a$ is a subset of $P$, and hence $a\in P$, that is
not the case). Denote by  $a'$ any point of  $A\setminus a$, which
does not belong to  $P$. Having expanded  $P\cup a$ to a maximal
convex set not containing $a'$, we obtain a piece $P'$ of $a'$,
which contains  $P\cup a$. $\blacksquare$\medskip

    {\bf  Definition.} The SS-rationalization  $(\mathcal P (f),\subset ,\phi
)$ is said to be the  {\em canonical rationalization} of a Plott
function $f$.\medskip

  {\bf  Example $1''$.} Let  $f$ be rationalizable by a partial
order $R$. Then the corresponding convex geometry $\mathcal F $
consists of ideals of the partial order. Since a lattice of ideals
$\mathcal F $ is closed under the union, for any element $x$ there
exists a maximal ideal $AF(x)$, which does not contain $x$. That is
any element $x$ has a unique piece $AF(x)$. We obtain that $\mathcal
P (f)$ coincides (as a poset) with $X$.\medskip

Now we are going to show that the SS-rationalization  $(\mathcal P
(f),\subset ,\phi )$ is not only canonical (that is uniquely
constructed by  $f$), but in some sense it occupies a ''central''
location among all SS-rationalizations of $f$. A perfect situation
would occur  if for any SS-rationalization $(Y,\le ,\psi )$ of $f$
would exist a mapping  $\alpha  \colon   Y \to \mathcal P (f)$,
such that $\phi \circ \alpha =\psi $ and $\alpha_* (g)=\tilde f$.
Unfortunately, this is not the case, as the following example from
(\cite{gree}) demonstrates.\medskip

{\bf  Example 6}.  Let  $X=\{a,b,c,d\}$ and let us define the Plott
function $f$ by the lattice of convex sets

 \unitlength=1mm
\special{em:linewidth 0.4pt} \linethickness{0.4pt}
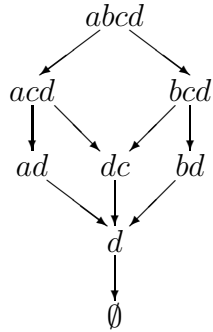
\begin{picture}(70.00,50)
\put(60.00,5.00){\makebox(0,0)[cc]{$\emptyset$}}
\put(60.00,15.00){\makebox(0,0)[cc]{$d$}}
\put(60.00,25.00){\makebox(0,0)[cc]{$dc$}}
\put(49.00,25.00){\makebox(0,0)[cc]{$ad$}}
\put(70.00,25.00){\makebox(0,0)[cc]{$bd$}}
\put(70.00,35.00){\makebox(0,0)[cc]{$bcd$}}
\put(49.00,35.00){\makebox(0,0)[cc]{$acd$}}
\put(60.00,45.00){\makebox(0,0)[cc]{$abcd$}}
\put(58.00,43.00){\vector(-4,-3){8.00}}
\put(62.00,43.00){\vector(4,-3){8.00}}
\put(68.00,33.00){\vector(-1,-1){6.00}}
\put(52.00,33.00){\vector(1,-1){6.00}}
\put(70.00,33.00){\vector(0,-1){6.00}}
\put(49.00,33.00){\vector(0,-1){6.00}}
\put(51.00,23.00){\vector(4,-3){8.00}}
\put(68.00,23.00){\vector(-1,-1){6.00}}
\put(60.00,13.00){\vector(0,-1){6.00}}
\put(60.00,23.00){\vector(0,-1){6.00}}
\end{picture}

In other words, the element $d$ is chosen only from  $\{d\}$. For
sets $A\neq X$ and $A\neq X\setminus d$,  the choice set is given
by $f(A)=A\setminus d$. Finally, $f(X)=f(X\setminus d)=\{a,b\}$.
The poset $\mathcal P (f) $ is the poset

 \unitlength=1mm \special{em:linewidth 0.4pt}
\linethickness{0.4pt}
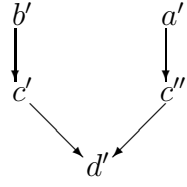
\begin{picture}(70.00,30.00)
\put(60.00,5.00){\makebox(0,0)[cc]{$d'$}}
\put(70.00,15.00){\makebox(0,0)[cc]{$c''$}}
\put(50.00,15.00){\makebox(0,0)[cc]{$c'$}}
\put(50.00,25.00){\makebox(0,0)[cc]{$b'$}}
\put(70.00,25.00){\makebox(0,0)[cc]{$a'$}}
\put(69.00,23.00){\vector(0,-1){7.00}}
\put(49.00,23.00){\vector(0,-1){7.00}}

\put(51.00,13.00){\vector(1,-1){7.00}}
\put(69.00,13.00){\vector(-1,-1){7.00}}
\end{picture}

However, if we consider a weaker poset  $Y$

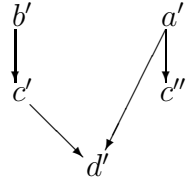
\unitlength=1.00mm \special{em:linewidth 0.4pt}
\linethickness{0.4pt}
\begin{picture}(70.00,30.00)
\put(60.00,5.00){\makebox(0,0)[cc]{$d'$}}
\put(70.00,15.00){\makebox(0,0)[cc]{$c''$}}
\put(50.00,15.00){\makebox(0,0)[cc]{$c'$}}
\put(50.00,25.00){\makebox(0,0)[cc]{$b'$}}
\put(70.00,25.00){\makebox(0,0)[cc]{$a'$}}
\put(49.00,23.00){\vector(0,-1){7.00}}

\put(51.00,13.00){\vector(1,-1){7.00}}

\put(69.00,23.00){\vector(-1,-2){8.00}}
\put(69.00,23.00){\vector(0,-1){7.00}}
\end{picture}

\noindent and the obvious mapping from $Y$ into  $X$ ($c'$ and $c''$
are sent into  $c$, $a'$ into $a$ and so on), then we obtain another
SS-rationalization of   $f$. The mapping  $\alpha $ is obviously the
identity mapping, but the corresponding Plott function $g$ (on $Y$)
is different from $\tilde f$ ($d'\in g(\{d',c''\})$ but $d'\not\in
\tilde f(\{d',c''\})$).

However, the following weaker property holds.\medskip

    {\bf  Proposition 3.} {\em  Let  $(Y,\le ,\psi )$ be an
SS-rationalization of a Plott function $f$. Then there exists a
mapping $\alpha \colon   Y \to \mathcal P (f)$, such that $\phi
\circ \alpha =\psi $ and $\alpha_*(g)\ge \tilde f$.}\medskip

Before proving this proposition we state its important corollary,
that the canonical rationalization is of minimal cardinality.
\medskip

    {\bf Corollary.} $|Y|\ge |\mathcal P (f)|$.\medskip

     In fact, because the support of  $\tilde f$ is the whole
$\mathcal P $, the support of  $\alpha_*(g)$ is also the whole
$\mathcal P $, and hence the mapping  $\alpha $ is a surjection.
$\blacksquare$\medskip

{\bf Proof of Proposition 3}. For a point  $y\in Y$ denote by
$AF(y)$ the set of points  $y'\in Y$ which are not dominated (with
respect to $\le $) by $y$. In other words, $AF(y)$ is the complement
in $Y$ to the principal filter $F(y)=\{y', y'\ge y\}$. The set
$AF(y)$ is an ideal with respect to the partial order $\le $ on $Y$.
Therefore, the set  $\psi_+ (AF(y))$ is a convex set in  $X$. This
set does not contain $\psi (y)$, therefore, there exists a piece $P$
of the element  $\psi (y)$, which contains  $\psi_+ (AF(y))$. We
define $\alpha  \colon   Y \to \mathcal P (f) $  by setting  $\alpha
(y)=P$ (if there are several pieces containing $\psi_+ (AF(y))$, we
pick up any). Because $\phi (P)=\psi (y)$, we obtain $\phi \circ
\alpha =\psi $.

In order to show  $\alpha_*(g)\ge f$ it suffices to check that the
direct image (under  $\alpha $) of the partial order $\le $ (on $Y$)
is weaker than the partial order  $\subset $ on $\mathcal P (f)$.
That is we have to check that if  $P$ and  $Q$ are two different
pieces and any element of  $\alpha ^{-1} (P)$ is dominated (with
respect to  $\le$) by some element of $\alpha^{-1}(Q)$, then
$P\subset Q$.

Let us prove this claim. Let $P$ be a piece of  $x$. Then $x$ is an
extreme point of  $P\cup x$, and hence  $x\in f(P\cup x)$. Since
$f=\psi_*  (g)$, we have, by the definition of the direct image, a
point $y$ of  $\psi^{-1}(x)$ which is not dominated by any point of
$\psi^{-1}(P)$. That is  $\psi^{-1} (P)\subset AF(y)$, or,
equivalently,  $P$ is a subset of the full image of $AF(y)$,
$P\subset \psi_+ (AF(y))$. Because the convex set $\psi_+ (AF(y))$
does not contain the point $\psi (y)=x$  and $P$ is a piece of $x$,
we have the equality $P=\psi_+ (AF(y))$ indeed. Thus, we have
$\alpha (y)=P$.

Now, any element of $\alpha ^{-1}(P)$ (and, hence, $y$) is dominated
by at least one element of $\alpha^{-1}(Q)$. Denote by $y'$ the
element of  $\alpha^{-1}(Q)$, which dominates $y$, that is $y<y'$
and $\alpha (y')=Q$. Due to the domination $y<y'$, we have the
inclusion  $AF(y)\subset AF(y')$, and, hence, $\psi_+(AF(y))\subset
\psi_+ (AF(y'))$ holds. By the definition of the mapping  $\alpha $,
$\alpha (y')=Q$ implies the inclusion $\psi_+(AF(y'))\subset Q$.
Since $P=\psi_+(AF(y))$, we obtain the desired inclusion $P\subset
Q$. $\blacksquare$

\section{Convex structure on {\bf SW}(X)}

In this section we define a convex structure on the set ${\bf
SW}(X)$ of simple words over $X$ (or on the set ${\bf LPF}(X)$ of
linear Plott functions).

Let $X$ and $Y$ be two disjoint sets, and let $f$ and $g$ be choice
functions on $X$ and $Y$ correspondingly. The {\em direct sum} of
$f$ and $g$ is the following choice function $f\coprod g$ on
the disjoint union $X\coprod Y$: for a subset $A\coprod B$ of $X\coprod Y$
$$
(f\coprod g)(A\coprod B)=f(A)\coprod g(B).
$$
Obviously, $f\coprod g$ is a Plott function if $f$ and $g$ are Plott
functions.

Let now $f$ and $g$ be linear Plott functions corresponding to
simple words $w$ and $v$ over $X$ and $Y$. We want to describe the
socle of $f\coprod g$. Without loss of generality we can assume that
$w$ and $v$ are  words of maximal length (that is the supports of $f$ and $g$
are equal to $X$ and $Y$, respectively). Let $\le_X$ and $\le_Y$ denote the
corresponding linear orders on $X$ and $Y$. Then we can consider the
partial order $R=\le_X \coprod \le_Y$ on $X\coprod Y$. As we have yet
seen in Example $1'$, the socle of $f\coprod g$ consists of all
linear extensions of the partial order $R$.\medskip

{\bf Definition.} A word (over the alphabet $X\coprod Y$)
corresponding to a linear extension of the partial order $R$ is
called a {\em shuffle} of $w$ and $v$. The set of shuffles of $w$
and $v$ will be denoted by $Sh(w,v)$.\medskip

A shuffle has a form
$$
w_1v_1w_2v_2\ldots w_1v_k,
$$
where $w=w_1w_2\ldots w_k$ and $v=v_1v_2\ldots v_k$; some of
sub-words $w_i$, $v_j$ are allowed to be empty. In words, a shuffled
word of two words is the following composition: To compose it we
have to take an initial piece of one of the words, then to append it
by an initial piece of  another word, then to return to the first
word and to take the next piece starting from the interruption
place, then again to switch to the another word and so on. We
illustrate this on the following example, let $w=xyz$ and $v$=abcd.
Then the word $x{\rm ab}y{\rm cd}z$ is one of the shuffles.\medskip

The shuffle operation enables us to define a (multi-valued) operation
$$
Sh: {\bf SW}(X)\times {\bf SW}(Y) \Longrightarrow {\bf
SW}(X\coprod Y).
$$
Obviously, it is commutative. But it is also associative. To see
this we note that the preceding construction can be done for
arbitrary number of words $w_i$ over disjoint sets $X_i$, $i\in
I$.\medskip

{\bf Lemma 4.} {\em Let $w$ be a shuffle of simple words
$w_1,w_2,...,w_n$ over $X_1,...,X_n$. Then there exists a shuffle
$w_{2...n}$ of $w_2,...,w_n$ such that $w$ is a shuffle of $w_{1}$
and $w_{2...n}$. }\medskip

{\bf Proof.} The linear order $\le$ corresponding to $w$ is an extension of
the order $\le_1\coprod ...\coprod \le_n$ on
$X_1\coprod...\coprod X_n$. Define $\le _{2...n}$ as the restriction
of $\le$ on the set $X_2 \coprod ...\coprod X_n$. Obviously, $\le$
is an extension of the disjoint union of $\le _{2...n}$ and $\le_1$. Now define
$w_{2...n}$ as the simple word corresponding to $\le _{2...n}$.
$\blacksquare$\medskip

Suppose now that $Y$ is a copy of $X$, and $\delta: X\coprod X \to
X$ is the co-diagonal mapping. A {\em melange} of two
simple words $w_1$ and $w_2$ over $X$ is a simple word
$\delta_{\sharp}(w)$ where $w$ is a shuffle of $w_1$ and $w_2$.

Let us give an example. Suppose we have two words $xyzab$ and
$zacyd$. Then the word $zaxycdb$ is a melange and it is  obtained by the
simplification of the following shuffle
$(za)(xyz)(cy)(a)(d)(b)$.
\medskip

Similarly we can define a melange of any family $w_i$, $i\in I$, of
simple words.\medskip

{\bf Proposition 4}. {\em Let $l_i$ be the linear Plott functions
corresponding to simple words $w_i$, $i\in I$. Then the set of
melanges of $w_i$ coincides with the socle of $\vee_{i\in
I}l_i$.}\medskip

{\bf Proof}. Let $\delta$ be the co-diagonal mapping of $X\times I$ onto
$X$, $\delta (x,i)=x$. It is obvious that $\delta_*(\coprod_i
l_i)=\vee_{i\in I}l_i$. Now the assertion follows from Theorem 1.
$\blacksquare$\medskip

As well as the shuffle, the melange is (multi-valued) commutative
and associative operation on the set ${\bf SW}(X)$. Using this
operation we can define a convex structure on ${\bf SW}(X)$.

Suppose that $w$ and $v$ are two simple words over $X$. Let
$co(w,v)$ be the set of all prefixes of all melanges of $w$ and
$v$. In other terms, $co(f,g)$ is the basement of $l_w\vee l_v$.
We may understand $co(w,v)$ as a ''segment'' joining the points
$w$ and $v$ in ${\bf SW}(X)$.\medskip

{\bf Definition.} A subset $C$ of ${\bf SW}(X)$ is said to be {\em
convex} if it contains $co(w,v)$ for every $w,v \in C$.\medskip

Lemma 4 has the following immediate consequence:\medskip

{\bf Corollary.} {\em Let $C$ be a convex subset of ${\bf SW}(X)$
and $c_1,...,c_n\in C$. If a simple word $w$ is a melange of
$c_1,...,c_n$ then $w\in C$.}\medskip

Indeed, by Lemma 3 $w$ is a melange of $c_1$ and $c_{2...n}$ where
$c_{2...n}$ is a melange of $c_2,...,c_n$. By induction $c_{2...n}$
is in $C$. Then $w$ is in $C$. $\blacksquare$\medskip

{\bf Theorem 4.} {\em Let $C\subset {\bf SW}(X)$. The following
two assertion are equivalent:

1) $C$ is a convex subset of ${\bf SW}(X)$;

2) $C=Bas(f)$ for some Plott function $f$ on $X$.}\medskip

{\bf Proof}. It is almost obvious that the basement $Bas(f)$ of
any Plott function is a convex subset of ${\bf SW}(X)$. Indeed,
let $w$ and $v$ be two simple words in $Bas(f)$, and let $l_w$ and
$l_v$ be the corresponding linear Plott functions. By the
definition this means that $l_w \le f$ and $l_v \le f$. Therefore
$l_w\vee l_v \le f$ and $co(w,v)=Bas(l_w\vee l_v) \subset Bas(f)$.

Conversely, let $C$ be a convex subset of ${\bf SW})X$. Let us
define $f$ to be the join of the linear Plott functions $l_c$
where $c\in C$. Obviously $C\subset Bas(f)$. Let us check the
inverse inclusion. Let $w$ be in $Bas(f)$ and even let $w$ be in
the socle of $f$. By Proposition 4, $w$ is a melange of several
words $c_1,...,c_n\in C$. By the previous Corollary $w\in C$.
$\blacksquare$\medskip

Thus, we obtain a bijection between the Plott functions on $X$ and
the convex subsets of the convex space  ${\bf SW}(X)$. This
bijection is compatible with the order structure. In particular
the lattice  ${\bf PF}(X)$ of Plott functions on  $X$ is
isomorphic to the lattice of convex sets of the set ${\bf
SW}(X)$.\medskip

{\bf Remark}. Because of this theorem and the bijection between
Plott functions and convex geometries, we obtain the bijection
between the lattice ${\bf CG}(X)$ and the lattice of convex
subsets of ${\bf SW}(X)$. Let us explain how to construct the
convex geometry $\mathcal F$ corresponding to a given convex
subset $C$ of ${\bf SW}(X)$. $\mathcal F$ consists of the
complements of the sets supp$(w)$ where $w$ runs over $C$ and
supp$(w)$ denotes the set of letters of a word $w$.\medskip

As a consequence of Theorem 2, we obtain a more explicit description
of the meet of a pair of  Plott functions $f_1$ and $f_2$. Let  $C_1$
and $C_2$ be the basements of $f_1$ and $f_2$. Then the function
$f_1 \wedge f_2$ corresponds to the convex set $C_1 \cap C_2$.
Similarly, $f_1 \vee f_2$ corresponds to the convex hull of $C_1
\cup C_2$.\medskip

{\bf   Example 7.} Let $X=\{a,b,c\}$. Consider two linear Plott
functions $f_1$ and $f_2$, which correspond to words $abc$ and
$bac$. The basement of $f_1$ consists of words: $abc$, $ab$, $a$ and
empty word $\emptyset$. The basement of $f_2$ consists of $bac$,
$ba$, $b$ and $\emptyset$. The intersection of these basements
consists of the single element $\emptyset$. That is $f_1 \wedge f_2=
{\bf 0}$.\medskip

    {\bf Example 8.} Let $X$ be as in the above example.
Consider other two orders  $abc$ and $acb$, and let $f_1$ and
$f_2$ be the corresponding linear Plott functions. The
intersection of their basements consists of two elements: the
empty word and the word $a$. Hence, $f_1 \wedge f_2$ is Plott
function with the support $\{a\}$ (the choice of a subset
$A\subset X$ is equal to $a$, if $a\in A$ and empty set
else).\medskip

{\bf Remarks}. 1. Similarly one can define a convex structure  on
the set ${\bf L}(X)$ of linear orders on $X$: a set $C \subset
{\bf L}(X)$ is convex if it contains all melanges of its elements.
We have used ${\bf SW}(X)$ in order to tame arbitrary Plott
functions.

2. The lattice of convex sets in  ${\bf SW}(X)$ is not a convex
geometry. (In other words, the lattice of all convex geometries is
not a convex geometry.) That might be shown on the three-elements
set $X$. Consider the whole set ${\bf SW}(X)$. It is a convex set.
But this set have no extreme points. In fact, the word $bca$ is a
melange of the words $bac$ and $cba$. Similarly, any (complete)
word is a melange of two other words. Therefore ${\bf
SW}(X)\setminus w$ is not convex for any word $w$.

Another example. Let us consider the segment $co(xzy, zxy)$. It
consists of seven words $xzy, \ xz, \ x, \emptyset, \ z,\ xx, \
zxy$. But the word $xzy$ is a melange of $x$ and $zxy$. Therefore it
is not an extreme point of this segment. Similarly the order $zxy$
is not an extreme point of this segment.

\section{The inverse image}

     Here we define the inverse image functor, the right conjugate
to the direct image; see a definition in \cite{G-Z}.

     Let $\phi :X \to Y$ be a mapping of sets and let $g$ be a
choice function on $Y$.\medskip

     {\bf Definition.} The {\em inverse image} $\phi^*(g)$ of $g$ is
the following Plott function on $X$:
$$
     \phi^* (g)=\wedge(f\in {\bf PF}(X), \phi_* (f)\le g).
                    $$

As the join of Plott functions, $\phi^*(g)$ is a Plott function.
Actually, $\phi^*(g)$ is the maximal Plott function on $X$ whose
direct image is $\le g$. In particular, if $\phi=id_X$ is identical
mapping of $X$, $\phi^*(g)=g^{\sharp}$ is the Plottization of $g$.

It is clear that
                 $$\phi_*  (\phi^*  (g))\le g.$$

     {\bf Example $4'''$.} Let $g={\bf 1}_T$ be the function of the
"identical" choice from a set $T\subset Y$. Then $\phi^*(g)={\bf
1}_S$ where $S=\phi^{-1}(T)$.

     Indeed, on one hand, $\phi_*({\bf 1}_S)={\bf 1}_T \le g$
so that $\phi^*(g) \ge {\bf 1}_S$. On other hand, the support of
$\phi^*(g)$ should be contained in $S$. Hence, $\phi^*(g) \le {\bf
1}_S$.\medskip

    {\bf  Proposition 5. (Conjugation)}. {\em  Suppose $\phi :X \to Y$
is a mapping, $f$ is a Plott function on $X$, and $g$ is a choice
function on $Y$. The following assertions are equivalent:

     1) $\phi_*  (f)\le g$,

     2) $f\le \phi^*  (g)$.}\medskip

     Proof. 1) implies 2) by the definition. Inversely, let $f\le
\phi^*(g)$. Then applying $\phi_*$ to the both sides, we obtain $\phi_*  (f)\le
\phi_* (\phi^*  (g))\le g$. $\blacksquare$\medskip

Suppose now that $g$ also is a Plott function. The inverse image
defines a mapping $\phi^* : {\bf PF}(Y) \to {\bf PF}(X)$ of sets
(and even a morphism of posets). It is rather obvious that the
inverse image is a (contravariant) functor from the category of
finite set to the category of posets.

As we know, there is a bijection between Plott functions and
convex sets in ${\bf SW}$. Which convex sets in ${\bf SW}(X)$
correspond to the inverse images? An answer is given by the
following
\medskip

    {\bf Proposition 6.} {\em Let $g$ be a choice function on $Y$.
    Then $Bas(\phi^*  (g))= \phi_{\sharp}^{-1} Bas(g)$.}\medskip

{\bf Proof}.   Suppose that $l$ is a linear Plott function on $X$. If
$l \le \phi^*  (g)$ then by Proposition 5, $\phi_*(l)\le g$.
Moreover $\phi_*(l)$ is a linear Plott function (Lemma 2). Therefore
it belongs to the basement of $g$. Inversely, if $\phi_*(l)\le g$
then by the definition $l\le \phi^* (g)$. $\blacksquare$\medskip

Since the meet of Plott functions corresponds to the intersection of
the basements, and $\phi_{\sharp}^{-1}$ commutes with the
intersection, we obtain the following\medskip

     {\bf Corollary.} {\em  The functor $\phi^*$ commutes with
$\wedge$.}$\blacksquare$\medskip

    {\bf Proposition 7. (The projection formula). } {\em Let $\phi
:X \to Y$ be  a map of sets, let $f$ be a Plott function on $X$,
and let $g$ be a Plott function on $Y$. Then}
$$
 \phi_*(f\wedge \phi^*(g))=\phi_* (f)\wedge g.
$$

{\bf Proof}. The inequality $\le$ follows from monotonicity of the
operator $\phi_*$. Let us prove the inverse inequality $\ge$.
Suppose that $l\in Bas(\phi_*(f)\wedge g)$. Then $l\in
Bas(\phi_*(f))$ and $l\in Bas(g)$. By Theorem 2, there exists a
linear Plott function $l'\le f$ such that $\phi_* (l')=l$. Since
$l\le g$, $l'\le \phi^*(g)$. Therefore $ l'\le f\wedge \phi^*(g)$
and $l=\phi_*(l')\le \phi_* (f\wedge \phi^*(g))$.
$\blacksquare$\medskip

     In particular, if $f={\bf 1}_X$ then we obtain the following
formula (where $g$ is an arbitrary Plott function on $Y$):
$$
       \phi_*(\phi^* (g))={\bf 1}_{\phi(X)}\wedge
g.
$$

     For example, if $\phi$ is a surjective map then $\phi(X)=Y$
and $\phi_*(\phi^*(g))=g$ for every Plott function $g$ on $Y$. Note,
that the inequality $\phi_* (\phi^*(g))\le g$ can be strict when
$\phi$ is not surjective.\medskip

    {\bf Example 9.} Let $Y$ consists of three elements $a,b,c$
and let $g$ be the linear Plott function corresponding to the word
$abc$. Let $X=\{b,c\}$ and let $\phi $ be the natural inclusion of
$X$ into $Y$. Because we cannot make a prefix of the word $abc$ from
letters $b$ and $c$, we obtain that the inverse image of $g$ is
equal to ${\bf 0}_X$. Therefore $\phi_*(\phi^* (g))={\bf
0}_Y$.\medskip

We give below two applications of the inverse image: to a
construction of direct products of Plott functions, and to a
construction of natural transformations of Plott functions.\medskip

{\em The direct product}

Let $X$ and $Y$ be sets, and let $\alpha$ and $\beta$ be the natural
projections of $X\times Y$ onto $X$ and $Y$. Let $f$ and $g$ be
Plott functions on $X$ and $Y$, correspondingly. Then we define the {\em
direct product} $f\prod g$ as the following Plott function on
$X\times Y$:
$$
f\prod g=\alpha^* (f) \wedge \beta^* (g).
$$

{\bf Example 10.} Let $X=\{x,x'\}$ and $Y=\{y,y'\}$. Suppose that
$f$ is the linear Plott function corresponding to the word $xx'$;
similarly $g$ corresponds to the word $yy'$. Let us compute the
direct product $f\prod g$. For this we describe its socle. The socle consists
A simple word $w$ (over the alphabet $X\times Y$) belongs to the socle iff its
projection onto $X$ is equal to $xx'$ and its projection onto $Y$ is
equal to $yy'$. That means that any such $w$ has  $(x,y)\in
X\times Y$ as its beginning letter. Therefore (for a set $Z\subset X\times Y$)
$$
(f\prod g)(Z)=\left\{\begin{array}{cl}
  \{(x,y)\} & \text{if} \ (x,y)\in Z,  \\
  Z & \text{if} \ (x,y)\notin Z.
\end{array}\right.
$$

{\em Plott correspondences}

Using the direct and inverse images and the meet
$\wedge$, one can construct very general transformations of Plott
functions. Namely, we call a {\em Plott correspondence} from $X$ to
$Y$ a Plott function $h$ on a set $Z$ and two maps $\phi:Z \to X$
and $\psi:Z \to Y$. If now $f$ is a Plott function on $X$ then
define
$$
                         h(f)=\psi_*  (h\wedge \phi^* (f)).
$$
As a result we obtain a map (of sets or of posets)
$$
h:{\bf PF}(X) \to {\bf PF}(Y).
$$

In the ''word'' interpretation this map looks as follows. We take
a convex subset $C$ in ${\bf SW}(X)$, lift it on ${\bf SW}(Z)$,
intersect with the basement of $h$, and then descent the
intersection onto ${\bf SW}(Y)$.

Note that we always can assume that $Z=X\times Y$. Indeed, let
$\alpha$ and $\beta$ be the natural projections of $X\times Y$ onto
$X$ and $Y$, and let $\pi:Z\to X\times Y$ be a mapping such that
$\phi=\alpha\circ \pi$ and $\psi=\beta\circ \pi$. Then
$$
h(f)=\psi_*  (h\wedge \phi^* (f))=\beta_*\pi_*  (h\wedge
\pi^*\alpha^* (f))=\beta_*(\pi_*(h)\wedge
\alpha^*(f))=(\pi_*(h))(f).
$$
That is, as the operators from {\bf PF}(X) to {\bf PF}(Y),
$h=\pi_*(h)$.

For example, if $Z$ is a subset of $X\times Y$, we can take $h={\bf
1}_Z$.\bigskip

{\bf Acknowledgments.} We want to thank A.Slinko for helpful
discussions and  advices on the subject.

\end{document}